\documentclass[12pt]{article}
\usepackage{amsfonts,amsmath}%,latexsym,amssymb}
\usepackage[round]{natbib}

\newcommand{\bZ}{{\bf Z}}
\newcommand{\bX}{{\bf X}}
\newcommand{\bx}{{\bf x}}
\newcommand{\bj}{{\bf j}}

\newcommand{\G}{{\cal G}}
\newcommand{\B}{{\cal B}}
\newcommand{\cP}{{\cal P}}
\newcommand{\cH}{{\cal H}}
\newcommand{\Hyp}{{\cal H}_{k}^\lambda}

\newcommand{\N}{\mathbb N}
\newcommand{\R}{\mathbb R}

\newtheorem{theorem}{Theorem}
\newtheorem{lemma}{Lemma}
\newtheorem{corollary}{Corollary}

\newtheorem{definition}[theorem]{Definition}
\newenvironment{proof}{Proof:}{$\Box$}

\date{February 7, 2003}

\title{Selection Criterion for  Log-Linear Models
Using Statistical Learning Theory}

\author{Daniel J. L. Herrmann \\
\small Max-Planck-Institut f\"{u}r biologische Kybernetik\\
\small Spemannstr. 38, D-72076 T\"{u}bingen, Germany \\
\small daniel.herrmann@tuebingen.mpg.de
\\ \\
Dominik Janzing \\
\small Institut f\"{u}r Algorithmen und Kognitive Systeme\\
\small Universit\"{a}t Karlsruhe, Am Fasanengarten 5 \\
\small D-76131 Karlsruhe,  Germany \\
\small janzing@ira.uka.de }

\begin{document}

\maketitle

\abstract{Log-linear models are a well-established method for describing
statistical dependencies among a set of $n$ random variables.
The observed frequencies of the $n$-tuples are explained
by a joint probability such that its  logarithm is a sum of functions, where 
each function
depends on as few variables as possible.
We obtain for this class a new model selection criterion using nonasymptotic
concepts of statistical
learning theory. We calculate the VC dimension for the class of $k$--factor 
log-linear models.
In this way we are not only able to select the  model with the appropriate
 complexity, but obtain
also statements on the reliability of the estimated probability distribution.
Furthermore we show that the selection of the best model among a set of models
with the same complexity can be written as a convex optimization problem.
}
%%%%%%%%%%%%%%%%%%%%%%%%%%%%%%%%%%%%%%%%%%%%%%%%%%%%%%%%%%%%%%%%%%%%%%%%%%%%%%%%
%%
%

{\bf Key words:} Markov network, AIC, non-asymptotic statistic, VC dimension
structural risc minimization, convex optimization
%
%%%%%%%%%%%%%%%%%%%%%%%%%%%%%%%%%%%%%%%%%%%%%%%%%%%%%%%%%%%%%%%%%%%%%%%%%%%%%%%%
%%
%
%
\section{INTRODUCTION}
Suppose a scientist is interested in the relation between $n$ features
described by the random variables $X_1,X_2,\dots,X_n$.
Observing  dependencies among these features is an important part of all
scientific disciplines. The dependencies may be deterministic, e.g.,
\[
X_n=f(X_1,\dots,X_{n-1})
\]
or  (in the generic case) statistical, i.e., the probability for the event
$X_1=x_1,X_2=x_2,\dots,X_n=x_n$, for short $\bX = \bx$, is not the product of 
the
probabilities of the $n$ events $X_j=x_j$.
All information about these dependencies is contained in the
joint probability distribution $P$ that assigns the probability
$P(\bx) = P(x_1,x_2,\dots,x_n)$ to each $n$-tuple.
Notice that for large $n$ no ``reasonable'' sample size
is sufficient to determine the probabilities of all
$P(x_1,x_2,\dots,x_n)$ with good reliability since the
size of $\Omega :=\times_j \Omega_j$ (where each $\Omega_j$ denotes the set of 
possible values $x_j$)
is exponentially large in $n$.
This shows that we should not try to learn the joint probability distribution,
we rather have to develop inference rules for learning the
statistical dependencies {\it in a weaker sense}.
The task is therefore to estimate some properties of the joint distribution
described by empirical data of a given size.

In the last years methods to judge scientific theories by their predicted  
models
on real world data sets  become  more and more important
(Hagenaars 1994, Ishii-Kuntz 1997, Pitt and Myung 2002)
and since the accessible
computational power is increasing strongly this trend will proceed in future.

It is a well-known fact, see e.g.
(Akaike 1973, Hansen and Yu 2001),
that criteria describing the
goodness of fit for the model to the data set only
are not enough to judge the scientific relevance of the model. Additionally one 
has to take
the complexity of the model into account. Usually the number of free
parameters are treated as a measure for the model complexity and are 
incooperated in
criteria like in Akaike's   Information Criterion (AIC)
(Akaike 1973),
or  Bayesian Information Criterion of Schwarz (1978).
Bayesian
Model Selection
(Kass and Raftery 1995)
and Minimum Description Length
(Rissanen 1996, Hansen and Yu 2001)
take also the functional
form of the model into account. However, these criteria do not give any 
statements
on the reliability of the estimation.

In this article we derive reliability statements from
statistical learning theory (Vapnik 1998) and show that the trade-off between
goodness of fit and model complexity can be treated with
structural risk minimization.

Let us describe roughly the link between the ideas presented here and the 
question of
minimal description length  in information theory, see
(Barron 1991, Rissanen 1996)
and references therein.
In the latter approach a risk functional is minimized over a set of
probability distributions where each  probability distribution is penalized
according to its description length. More explicitly,  the risk functional 
consists
of two terms: the  ``empirical risk'' evaluates the ``goodness of fit''
and a regularization term penalizes the code length.
Barron and Cover (1991)
obtain bounds on the convergence
rate of the risk functional in the Hellinger distance.
The difference to our approach lies in the fact  the regularization function
can in principle be chosen arbitrarily, only  Kraft's inequality most hold;
whereas we use a penalty term coming directly
from the complexity of the considered class of log-linear models
and is based on statistical learning theory.
Vapnik and Chervonenkis (1971)
addressed the question of
model selection from empirical data in a probabilistic setting
(see also Vapnik (1999)
and references therein).
They introduced the so--called VC dimension of a function set and showed that
the VC dimension of the function set which the model can implement is a crucial 
quantity to
describe its complexity. In general the VC
dimension does not agree with the usual dimension or the number of free 
parameters.

From  the statistical learning  point of  view, a scientific theory delivers a 
prior on the class of models
under consideration. Then the model is chosen according to  some criterion 
measuring the fit on the
empirical data.
In this work we consider the class of all log-linear models and prior models 
with few statistical
dependencies among the random variables, i.e. we prefer less complex models in 
terms of
statistical independency as long as they describe the data well enough.
With this assumption we estimate the model (complexity) directly from the 
empirical
data rather than testing  one model against another  using test methods like
the Pearson or likelihood ratio chi-square
criterion, see e.g.
Christensen 1997, Goodman 1978.
Akaike (1973)
already formulated this view point, however he did not have the
conceptual tools of learning theory
(Vapnik 1995).

The article is organized as follows.
In Section \ref{sec:2} we introduce log-linear models and explain the
prior for model selection
leading to a risk functional on the data. We discuss the commonly used model 
selection criteria and
explain their shortcomings.
We present the idea of Markov networks and relate them to log-linear models.
In Section \ref{sec:Risc} we calculate the VC dimension of the models under 
consideration
and give the central result about the estimation of the risk functional for the 
true probability
distribution.
In Section \ref{sec:Inf} we use the idea of the structural risk minimization to
formulate the model estimation for log-linear models.
In Section \ref{sec:convex} we show that the optimal model is the solution of
a convex optimization problem.

%%%%%%%%%%%%%%%%%%%%%%%%%%%%%%%%%%%%%%%%%%%%%%%%%%%%%%%%%%%%%%%%%%%%%%%%%%%%%%%
%
%
%\section{THE IDEA OF LOG-LINEAR MODELS}
\section{LOG-LINEAR MODELS}
\label{sec:2}
\subsection{The Empirical Risk Functional}% for Log-linear Models}
In the work presented here we restrict our
attention to discrete random variables and assume that $\bX$ takes
values in a finite set $\Omega = \times_{j=1}^n\Omega_j$.
Assume a scientist tells us that he has found that the true joint
probability distribution $P$ is given by $P_t$ and we
would like to test whether he is right or not.
Suppose that his model is given by
the distribution that all $n$-tuples in $\Omega$ occur with equal probability.
Assume, in contrast,  that the {\it true}
distribution $P$ assigns the probability
$1/k$ to $k$ specific $n$-tuples ${\bx}_1,{\bx}_2,\dots,{\bx}_k
\in \Omega$ and that this $k$ $n$-tuples are chosen without any simple law.
Assume furthermore that $k$ is large compared to the available sample size
but small compared to the size of $\Omega$.
Then we have only little chance to recognize that his model is wrong.
Only if the $k$ $n$-tuples are selected by a simple law we would mistrust
the scientist's hypothesis.
This example illustrates that we need a quality criterion for models
which can be tested on a reasonable sample size.  We suggest to use following 
criterion.
For each $n$ tuple ${\bx}$ in our data set we give penalty points depending on 
how
likely its appearance was according to the hypothetical joint distribution.
A good choice, for instance, is the negative logarithm of the hypothetical
probability $P_t({\bx})$. Let $T$ be the training set, then we obtain
\[
-\frac{1}{\# T}\sum_{{\bx} \in T} \ln P_t({\bx})
\]
as the ``penalty'' for the model $P_t({\bx})$. Notice that this sum converges
in the large sample size limit to
\begin{equation}
  \label{eq:Risc}
R(P_t):=-\sum_{{\bx} \in \Omega} P({\bx})
\ln P_t ({\bx})\,.
\end{equation}
Note that
$R(P_t)$ is closely related to the Kullback-Leibler distance (Cover and
Thomas 1991)
and can never be smaller than the entropy
\[
S(P):=-\sum_{{\bx} \in \Omega} P({\bx}) \ln P({\bx})\,.
\]
Hence, even if the scientist has predicted the true measure exactly,
he will never get less than $S(P)$ penalty points.
If $S(P)$ is high, he is not punished for building {\it wrong}
models but for selecting the \textit{wrong} features, i.e.~features
with too few dependencies.
If the scientist considers additional random variables
(``features influencing the considered ones'')
more specific statements might be possible.
Reasonable research consists not only in observing
dependencies among a given set of features, it consists also
in {\it  selecting} statistically relevant features
yielding a low entropy of the joint distribution.
Therefore the ``empirical risk functional''
\begin{equation}
  \label{eq:empRisc}
  R_{emp}(P_t) := -\frac{1}{\# T}\sum_{{\bx} \in T} \ln P_t({\bx})
\end{equation}
is a good measure to {\it test} the quality of the distribution $P_t$
for large sample size. This risk functional is well-known for estimating 
probabilities
(Vapnik 1998).

If the space of all probability distributions is too large compared to
the sample size, we should not minimize the empirical risk over all
joint distributions. In order to infer a probability distribution from the
statistical data we give the space a class structure.
We use log-linear models (Christensen 1997, Goodman 1978)
to define suitable classes.
The simplest class of log-linear models is given in the case
when all random variables are statistically independent. Then the logarithm of 
the
joint probability distribution can be written as a sum of functions, each 
function
depends on one variable only
\[
\ln P_t(x_1,x_2,\dots,x_n)=\sum^n_{j = n} f_j(x_j)\,,
\]
with $f_j(x_j)= \ln P(x_j)$.
The next higher class is given by allowing two-variable interaction
terms
\[
\ln P_t(x_1,\dots,x_n)= \sum_{j\leq n} f_j(x_j) +\sum_{i,j \leq n}
g_{i,j}(x_i,x_j)\,.
\]
The next higher class contains terms with three-variable interactions, etc..
Probability distributions with the property that
their logarithm can be expressed as sum of functions with $k$ variables
are called $k$-factor distributions. In other words the data is described
by a \textit{$k$--factor log-linear model}.
The idea of log-linear models come from the fact that some variables
influence each other {\it directly}, expressed by the interaction terms,
whereas the other dependencies are caused {\it indirectly} by influencing
intermediate variables.
The hope that real world data are well explained by log-linear models
can be backed up by the idea that each variable is influenced only by a few
other variables directly. The graph of the variables with edges representing the
direct influences is simple. This idea can be formulated using {\it Markov 
networks},
see Section \ref{sec:Markov}.
\subsection{Model Selection and AIC}
\label{sec:akaike}
A simple way of choosing  the degree of the interaction terms in a log-linear 
model
goes as follows. We begin with the simplest model class and decide on the basis
of some significance test whether the data suggest to reject the model or not.
If so then we include two-variable interactions. The same significance test
is applied to the new model. The procedure ends if frequencies observed in the
data and the probabilities given by the model do not differ significantly  
according to
the significance test. Let us describe two popular significance tests
(Christensen 1997, Goodman 1978).

Let $j \in \{1,\dots,k\}$ be the probability space consisting of $k$
possible events. Let $p_1,\dots,p_k$ be the  probabilities of the model
and $m_1,\dots,m_k$ the observed frequencies.
The sum $l:=\sum_j m_j$ is the sample size. Then we calculate

\begin{enumerate}

\item {\bf Pearson chi-square Test:}
\[
X^2:=\sum_j \frac{(m_j - lp_j)^2}{lp_j}
\]

\item {\bf Deviance Test:}
\[
G^2:=2\sum_j m_j \ln \frac{ m_j}{lp_j}\,.
\]

\end{enumerate}
\textit{Remark:} $G^2/(2l)$ is the Kullback-Leibler distance (Cover and Thomas 
1991)
between the probability distribution $(p_j)_{j\leq k}$ and the relative 
frequencies
$(m_j/l)_{j\leq k}$. Furthermore the minimizer of $G^2$ and  $R_{emp}$ are the 
same, i.e.,
$G^2/(2l) - R_{emp}$ does not depend on  $(p_j)_{j\leq k}$.
\\

Under the assumption that  $(p_j)_{j\leq k}$ is the true probability 
distribution the values
of the random variables $X^2$ and $G^2$ are in the large sample limit
$\chi^2$-distributed (Christensen 1997).
Using this argument the model is rejected if
the $X^2$ or $G^2$ values are too large. The range of acceptable values for
$X^2$ and $G^2$ is based  on the following rule:

First we have to calculate the number of
degrees of freedom ($df$) of the model. ($df$) is given by
the number of possible events minus $1$ minus the number of
free parameters in the model class.
Given $df$ one can look up in a table with $\chi^2$ distribution whether
the probability that $X^2$ or $G^2$, respectively, is outside
a certain confidence interval.

However, the idea of rejecting a simple model only if the data contradicts the 
externally
given significance level leads to conservative models.
Suppose the outcome of the chi-squared test of a model yields
values of $\chi^2$ or $G^2$ within a certain confidence interval.
Nevertheless there might be
a model with only one interaction term more such that the $\chi^2$
or $G^2$ value is considerably decreased.
In this case we will certainly prefer the model that is {\it
a  little bit} more complex.
This example illustrates that the externally given significance level implicitly
determines the model selection.
Therefore it is crucial to find a more systematic way
for model selection and prevent in this way over- and underfitting.

Let $(P_\theta)_{\theta \in \R^k}$ be a set of probability distributions
with $k$ parameters and assume that the true probability distribution
$P$ lies in $(P_\theta)_{\theta \in \R^k}$.
The true probability distribution $P$ is the minimizer of $R$  and
for a sample size which is large compared to the size of $(P_\theta)_{\theta \in 
\R^k}$,
the empirical risk functional $R_{emp}$
is a good approximation for $R$. In this case we can search for the minimizer of
$R_{emp}$ instead. This leads to the Deviance test above.
However we are interested in the case where we cannot assure that the data set
is sufficiently large.
Then the minimizer of $R_{emp}$ is in general not a good estimation for the
minimizer of $R$.

Akaike (1973)
suggested to modify the empirical risk functional according
to the  number of free parameters of the model under consideration.
He showed under some regularity conditions on the mapping $\theta \mapsto 
P_\theta$, that
\[
R_{emp} + df/l
\]
is a consistent estimator of $R$, i.e., it has the same expectation value
as $R$. However minimizing $R_{emp} + df/l$ for a finite data
set give us no information about the true risk $R$.
In  Section \ref{sec:Risc} we obtain some bounds on the
true risk and suggest to  minimize this {\it guaranteed} risk rather than
the {\it expected} risk, see Section \ref{sec:Risc}.
\subsection{Why Log-linear models}
\label{sec:Markov}
In Section \ref{sec:2}  we have introduced log-linear model classes
without giving any justification for it.
In this section we explain why one should expect that
probability distributions generating real world data can be described by
$k$-factor models with small  $k$.
The key idea is that the distribution of a random variable $X_j$
is only influenced by the values a few others variables.
This kind of statistical relations can be encoded by a Markov network.

The presentation of Markov networks in this work follows Pearl (1988).
For each random variable $X_i$ we define the Markov boundary
$\B_i \subset \{X_1,\dots,X_{i-1},X_{i+1},\ldots,X_n\}$ as the smallest
set fulfilling
\begin{equation}
  \label{eq:Markov}
P(x_i |x_1,\dots,x_{i-1},x_{i+1},\dots,x_{n}) = P(x_i|\B_i )\,.
\end{equation}

Observe that this definition does not take any (temporal or causal) order of
the random variables into account.
It will therefore lead
to an undirected graph. This is in contrast to
the directed graphs in the context of
Bayesian networks (Pearl 1985) which are  for instance
useful to formalize
causal structure (Pearl 2000).

In case of strict positive probability distributions the Markov boundary is 
unique.
Actually, strict positivity of the probability distributions is also
a crucial
condition for our approach to learn probabilities
as we will see in Section \ref{sec:Risc}.
We now define the Markov network as the graph $\G$ with $n$ nodes also denoted 
by
$X_1,\dots,X_n$ and undirected edges between $X_i$ and $X_j$ if and only if
$X_j \in B_i$ (or equivalently $X_i \in B_j$).
In  other words, we can generate the Markov network for a distribution $P$ if
we connect each of the $n$ nodes with the element of its Markov boundary.
As a consequence we obtain that two random variables $X_i$ and  $X_j$, $i\not=j$ 
are
conditionally independent  with respect to a set
$\bZ \subset {\bX} \setminus \{X_i,X_j\}$, i.e.
\[
P(X_i | X_j, Z) = P(X_i|Z)\,,
\]
if and only if every path from $X_i$ to $X_j$
has one node in $\bZ$. In other words, if  $\bZ$ \textit{blocks} every path 
between
$X_i$ and $X_j$ then $X_i$ and $X_j$ are conditionally independent (with respect 
to $\bZ$).
A \textit{clique} in the graph $\G$ is a subgraph in which all nodes
are connected to each other. We have an (partial) order on the cliques given by 
the inclusion as
sets.

Let $\G$ be a Markov graph for a strictly positive probability distribution $P$ 
and
let $C$ be the set of maximal cliques in $\G$. Then one can show that $P$
factorizes as follows
\begin{equation}
  \label{eq:factor}
P(\bx)= \frac 1N \prod_{c \in C}  \phi_c(\bx_c)\,,
\end{equation}
where $\bX_c \subset {\bX}$ is the set of random variables in the cliques $c$.
$N$ is the normalization constant and $\phi_c$ are the \textit{compatibility 
functions}
w.r.t. the clique $c$, also called the potential functions.

Suppose two distributions have the same Markov graph, then they describe the 
same
(in)dependency  in  $\bX$.  In Section \ref{sec:Risc}  we will see that the
maximal number in a clique $\deg(\G) = \max_{c \in C} \# c$ are a
appropriate parameters
to subdivide the set of all distributions $\cP(\bX)$ into
\textit{ small} model classes.

Consider the case that every  clique has less than $k+1$ elements.
Then the distribution can be described by a $k$-factor model due to
the factorization in eq. (\ref{eq:factor}).
The Markov boundary can be understood in a rather literary sense if
the variables represent some quantities that are measured at different
positions in the real space. Then the belief to obtain
probability distributions that correspond to simple Markov networks
stems from the locality  principle: distant variables
influence each other only indirectly via other variables.

In statistical physics log-linear models can be justified even more directly.
Assume the variable $X_j$ describes the physical state
of particle $j$ and that the total energy of the system in the state
${\bf x}=(x_1,\dots,x_n)$
is a sum of functions depending
on at most $k$ particles, i.e.,
\[
E({\bf x})= \sum_{\bf j} f_{\bf j}(x_{j_1},\dots,x_{j_k})\,,
\]
where ${\bf j}$ runs over all $k$-subsets in $\{1,\dots,n\}$.
Then the thermodynamic equilibrium state, the so-called Gibbs distribution,
is up to a normalization factor directly given by
the exponential of the negative energy function.
Hence the equilibrium distribution is exactly described by a $k$-factor model.
In statistical physics one usually considers  two particle interactions only, 
therefore
$2$-factor models describe the equilibrium state.

Now we have justified why the class structure of log-linear models is natural.
However, we have to refine our hierarchy in some respect.
So far, we did not formulate any restrictions
on the range of the probability distribution of the models.
Intuitively it is clear that we can assign a small probability to an event
$\bx \in \Omega$ only if we have a very large data set. This means
that reasonable model selection should allow small probabilities
only for large data sets.
However,  the empirical risk functional $R_{emp}$ penalizes for small 
probabilities  only
if they are contained in the data set.
Therefore the value of $R_{emp}$ depends strongly on the specific data set if 
the probability
distribution is not bounded from below by a strict positive constant.
Actually, for deriving a bound for the risk functional $R$ in Section 
\ref{sec:Risc}
we precisely need to bound the probability in this way.
\begin{definition}
\label{def:family}
Let  $\cP_\lambda$ be the set of all  strict positive functions greater than 
$\lambda>0$.
Let $\cP^k$ be the set of probability distributions $P$ corresponding to a
$k$-factor model, i.e.,
\[
P(x_1,x_2,\dots,x_n)=\prod_{{\bf j}}
q_{{\bf j}}(x_{j_1},x_{j_2},\dots,x_{j_s})\,,
\]
where each  $q_{{\bf j}}$ is a positive function and
${\bf j}:=(j_1,j_2,\dots,j_s)$ runs over all
possible subsets of $\{1,2,\dots,n\}$ with $s$ elements.
Then we define the model class of degree $k$ by
\begin{equation}
  \label{eq:hypSpace}
   \cH_{k}^\lambda  = \cP_\lambda(\bX) \cap \cP^k\,.
\end{equation}
\end{definition}
\textit{Remark}: The values $q_{{\bf j}}(x_{j_1},x_{j_2},\dots,x_{j_s})$
cannot be interpreted as probabilities since they can be greater than $1$.

%
% However
%it also possible to use other norms and obtain bounds for the VC depending on 
the
%moment of   $\ln P(\bx)$.
%%%%%%%%%%%%%%%%%%%%%%%%%%%%%%%%%%%%%%%%%%%%%%%%%%%%%%%%%%%%%%%%%%%%%%%%%%%%%%%%
%
%
%
\section{LEARNING LOG-LINEAR \\ MODELS}
It is a well-known problem in any mathematical theory of learning
that a model can explain the training data well but does not fit the
data observed in the future if the model class is too large.
If the model class is small enough, one knows with high confidence that the
model which minimizes the error on the training data will also explain
the future data almost optimal within this model set.
\subsection{Risk estimation}
\label{sec:Risc}
For a set of real-valued  functions on an arbitrary set $\Omega$ a relevant
measure for the size of this set is  the so-called  VC  dimension
(Vapnik-Chervonenkis). First, let us define the VC dimension (Vapnik 1998)
of a set of indicator functions $(f_\alpha)$, i.e. $f_\alpha: \Omega \rightarrow 
\{0,1\}$.

\begin{definition}
Let $\Lambda$ be an index set of arbitrary cardinality.
Let  $(f_\alpha)_{\alpha \in \Lambda}$ be a set of indicator functions
on $\Omega$. Then the $VC$ dimension of $(f_\alpha)_{\alpha \in \Lambda}$
is the largest number $h$ such that there exist
$h$ points $x_1,x_2,\dots,x_h \in \Omega$ with the property
that for every function\\
$\chi: \{x_1,\dots,x_h\} \rightarrow \{0,1\}$
there exists a function $f_\alpha$ such that its restriction
to $\{x_1,x_2,\dots,x_h\}$ coincides with $\chi$.
\end{definition}
The VC dimension for a set of real-valued functions is defined  by
the VC dimension for the set of corresponding indicator functions.

\begin{definition}
Let $(f_\alpha)_{\alpha\in \Lambda}$ be a family of real-valued functions
on a set $\Omega$. Then the VC dimension of $(f_\alpha)_{\alpha\in \Lambda}$
is the VC dimension of the family of the indicator functions
$(\chi_\mu \circ f_\alpha)_{\mu \in \R, \alpha \in \Lambda}$, where
$\chi_\mu$ is $\chi_\mu(x)=0$ for $x< \mu$ and $\chi_\mu (x)=1$ for $x\geq \mu$.
\end{definition}

The notion of VC dimension plays a crucial role in statistical learning theory.
It is a measure which indicates whether a learning machines overfits or not.
Roughly speaking, if a family of functions which can be implemented by a 
learning machine
has small VC dimension then every function which fits well the training data
will fit with high probability the test data as well. The
following theorem is an Corollary  of Theorem 5.1 pp.192 in (Vapnik 1998).

\begin{theorem}
\label{theo:Lern}
Let $(f_\alpha)_{\alpha \in \Lambda}$ be a measurable set of bounded
real-valued functions on $\Omega$, $A \le f_\alpha(z) \le B$, and let the set of 
indicator functions have finite VC
dimension $h$. Let $\mu$ be a probability distribution on $\Omega \times \R$ and
let the data  $\bx_1,\bx_2,\dots,\bx_l$ be drawn according to  $\mu$
(independent and identically distributed).
Then we have with probability at least $1-\eta$
\[
\tilde{R}(f_\alpha) \le \tilde{R}_{emp} (f_\alpha) +  (B-A)
\sqrt{ \frac{h -\ln h + \ln 16l   - \ln \eta}{l} } \,,
\]
where $\tilde{R}_{emp} (f_\alpha) := \frac{1}{l}\sum_{j\leq l} f_\alpha
({\bf x_j})$
is the empirical risk and
$\tilde{R}(f_\alpha) :=\int f_\alpha ({\bf x}) d\mu({\bf x})$ the true risk.
\end{theorem}
\textit{Remark:} We have  $R(P_\alpha) = \tilde{R}(\ln P_\alpha)$ and
$R_{emp} (P_\alpha) = \tilde{R}_{emp} (\ln P_\alpha)$ when $(P_\alpha)$ is a set 
of
probability distributions.
Therefore we can use this theorem to calculate the generalization error of 
equation
(\ref{eq:Risc})  for the family of probability distributions defined in equation
(\ref{eq:hypSpace}).

First, let us interpret the generalization error in this context.
Suppose we have found a probability distribution $P_t$ in  some model class
that minimizes $R$ then we know that in this space $P_t$ is the closest 
distribution
to the true distribution in the Kullback-Leibler distance. Furthermore, we
can conclude with high confidence that future data will fit to this
distribution as well.

\iffalse
\begin{corollary}
Let $(P_\alpha)_{\alpha \in \Lambda}$ be a family of
probability functions
on $\Omega$ with VC dimension $h$. Then the probability that
for any $\alpha$
after $l$ runs
the empirical risk $R_{emp}(P_\alpha)$
deviates from the true risk
$
R(P_\alpha)
$
by more than
\[
\phi (h)  =2\sqrt{\frac{h(\ln(2l/h) +1)-\ln (\eta/4)}{l}}
\]
is less than $1-\eta$ for arbitrary $\eta >0$.
\end{corollary}

Intuitively, this corollary states the following.
Assume we have considered a family of probability distributions with small
VC dimension. Once we have found a distribution that fits well to the data
in the sense that $n$-tuples that has been predicted to be unlikely
have not occurred very often (i.e. the empirical risk is small), then
we can conclude with high confidence that future data will
fit to this distribution as well.

Based on VC dimension, one can use
the so-called {\it structural risk minimization principle}.
The idea is to define an increasing net of families of functions
such that $R_{emp} +\phi(h)$ is minimal.
One the one hand, a richer family of functions may allow
to reduce the empirical risk $R_{emp}$.
On the other hand, a richer family may have a higher VC dimension
and will therefore increase the confidence term $\phi(h)$.
The structural risk minimization principle helps to allow
families that are large enough to explain data well and small enough
to avoid overfitting and destroy the generalization ability.
\fi

Next we estimate the VC dimension of product distributions
since this will be the leading intuition for estimating the
VC dimension of the model class $\cH_{k}^\lambda$.

\begin{lemma}
Let  $X_1,X_2,\dots,X_n$ be discrete random variables where $X_i$ takes $m_i$ 
different
values in $\Omega_i$ for $i=1,\ldots,n$.
Then the VC dimension of the family of all product distributions
$\cP^1 = \{P(x_1,\dots,x_n)=P(x_1)\cdot\ldots\cdot P(x_n):\; P \text{ is a 
probability distribution} \}$
is given by  $N = 1 + \sum_{i=1}^n m_i -n$.
\end{lemma}
\textit{Remark}:
For the families $(\ln P)_{P\in \cP^1}$  and  $(P)_{P\in \cP^1}$ the set of 
indicator
functions coincides. Hence $(\ln P)_{P\in \cP^1}$  and  $(P)_{P\in \cP^1}$
have the
same VC dimension.

\vspace{0.3cm}

\begin{proof}
Let $x_{j;0}$ be an arbitrary element of the set $\Omega_j$.
Due to $P(\bx)=P(x_1)\dots P(x_n)$ we can
write the logarithm of the joint probability as
\[
\ln P(\bx)=\sum_j \ln (P(x_j)/P(x_{j;0})) +\sum_j \ln P(x_{j;0})
\]
We can characterize an $n$-tuple $\bx =(x_1,x_2,\dots,x_n)$  uni\-quely
by an $N=\sum_j (m_j-1)$ dimensional
 vector with entries $0$ and $1$ as follows:
The vector consists of $n$ blocks of dimension $(m_j-1)$ characterizing
$x_j$. Let $x_{j;0},x_{j:1},\dots,x_{j;m_j-1}$ be the elements of
$\Omega_j$ in an arbitrary ordering.
For  $x_j=x_{j;i}$ with $i\neq 0$
let the $i$-th entry of the $j$-th block be $1$ and
the other entries of the $j$-th block be zero.  For $x_j=x_{j;0}$
let all entries of the $j$-th block be zero.

For each $n$-tuple the logarithm of its probability is  given as follows.
Let $c^{\bx} \in \R^N$ be the vector corresponding to  ${\bx}
:=(x_1,x_2,\dots,x_n)$.
Then we have
\[
\ln P({\bx})=\sum_l c^{\bx}_l f_l + d =\langle c^{{\bx}}, f\rangle +d
\]
where the coordinates $f_l$ of $f$ are
given by the values $\ln (P(x_{j;i})/P(x_{j;0}))$ with $j\leq n$ and
$1\leq i\leq m_{j}-1$.
The constant term $d$ is defined as
$d:=\sum_j \ln P(x_{j;0})$.
This shows that $\ln P({\bx})$
can be written as an affine functional in $\R^N$.
Hence the VC dimension of $\cP$ is at most $N+1$
(see (Vapnik 1998), Example in Chapter 5.2.3).
To show that it is not smaller than $N+1$ note that
there is no restriction at all to the set of all possible vectors $f$
if $d$ is chosen appropriately. This can be seen as follows.
Let $g_1,g_2,\dots,g_{m_j-1}$ be  the entries of $f$ in the $j$-th block.
Then set
\[
P(x_{j;i}):= \frac{\exp(g_i)}{\sum_{1\leq i \leq m_j-1}\exp( g_i)+1}
\]
and
\[
P(x_{j;0}):=\frac{1}{\sum_{1\leq i \leq m-1}\exp( g_i)+1}\,.
\]
Now we have to show that
one can find $N+1$ vectors that can be classified
in all $2^{N+1}$ possible ways such that the vectors correspond to
$n$-tuples.
Consider the vectors $e_j$  having $1$ at the position $j$ and $0$ elsewhere.
Obviously, each vector $e_j$ corresponds to a possible $n$-tuple.
As $N+1$-th vector we choose the origin $(0,0,\dots,0)$. It corresponds
to the $n$-tuple $(x_{1;0},x_{2;0},\dots,x_{n;0})$.
These vectors  can be classified in all $2^{N+1}$ possibilities
as follows. Choose an arbitrary indicator function $\chi$ on these
$N+1$ points
and set $f:=(\pm 1,\pm 1,\dots,\pm 1)$
with positive sign at position $j$ if and only if $\chi(e_j)=1$.
Define $d$ as above
in order to ensure that
\[
{\bx}\mapsto \langle c^{{\bx}}, f\rangle +d
\]
defines the logarithm of a probability function.
Set $a:= d \pm 1/2$ and choose the sign positive  if and only if
$\chi ( (0,0,\dots,0))=1$.
Then the indicator function $\chi_a \circ w$ with
$w ({\bx}):= \langle c^{{\bx}} ,f \rangle +d $ coincides
with the desired classifier $\chi$.
Hence the VC dimension of product distributions is $N+1$.
\end{proof}

Now we  give a bound for the VC dimension of the model classes introduced in 
Section
\ref{sec:Markov}.
\begin{lemma}
\label{VCdim}
The VC dimension of $\cP^k$ is at most
\[
h_k:=\sum_{{\bf j}} m_{j_1} m_{j_2} \dots m_{j_k}\,,
\]
where $m_i$ is the size of the set $\Omega_i$.
\end{lemma}
\textit{Remark:} $h_s$ is at most $O(n^k)$
if $m_i$ is uniformly bounded.

\begin{proof}
Let  $P\in \cP^k$. Since $P$ factorizes we can characterize
the expression $\ln P$ uniquely by a
vector $f\in \R^h$ as follows: the vector $f$ consists of
${n \choose k}$ blocks where the blocks are indexed by $\bj$ given in the Lemma.
The block $\bj$ has size $ \prod_i m_{j_i}$ and the entries within the block 
have  the values
\[
\ln q_\bj (x_{j_1},x_{j_2},\dots,x_{j_k})\,,
\]
where $(x_{j_1},x_{j_2},\dots,x_{j_k})$ runs over all elements
of $\Omega_{j_1}\times \Omega_{j_2} \times \cdots \times \Omega_{j_k}$.
Each $n$-tuple ${\bx} =(x_1,x_2,\dots,x_n)$
can be characterized by a vector $c^{{\bx}}\in\R^h$ as follows:
In each block ${\bf j}$ there is exact one coordinate that corresponds
to the $s$-tuple $(x_{j_1},x_{j_2},\dots,x_{j_k})$. Set this entry to $1$
and the other entries of the block to $0$.
Then $\ln P({\bx})$ is given by the inner product of the vectors
$f$ and $c^{{\bx}}$.
We obtain the lemma since the VC dimension of the set of linear functionals
on $\R^h$ is $h$, see (Vapnik 1998).
\end{proof}

Obviously, the model class $\Hyp$ is a subset of $\cP^k$ and so we can
give an upper bound the VC dimension of $\Hyp$ with $h_k$ which
is independent of $\lambda$.
Therefore  Lemma \ref{VCdim} and Theorem \ref{theo:Lern} yields to the
following theorem.
\begin{theorem}
\label{theo:Lern2}
Let $\Hyp$ be the model class of degree $k$ and $h_k$
be the upper bound on
the VC dimension of $\Hyp$ given in Lemma \ref{VCdim}.
Then for a training data $\bx_1,\bx_2,\dots,\bx_l$ we have for $P_\alpha \in 
\Hyp$
with probability at least $1-\eta$
\begin{equation}
  \label{eq:risc}
  R(P_\alpha) \le R_{emp} (P_\alpha) + \phi (k,\lambda,\eta)\,,
\end{equation}
with $\phi (k,\lambda,\eta) := - \ln \lambda \;
            \sqrt{ \frac{h_k -\ln h_k + \ln 16 + \ln l - \ln \eta}{l} }$.
\end{theorem}
The task to find the  minimizer of  $R_{emp}$ in $\Hyp$ even for fixed
$\lambda$ and $k$ is computationally expensive. We will not treat this question 
here.
Theorem \ref{theo:Lern2} tells us how to estimate the risk $R(P_\alpha)$ for
a a priori chosen model class $\Hyp$. However we do not know which class to 
select.
This question is  addressed in the next section.
\subsection{Estimating Log-linear Models  from empirical data}
\label{sec:Inf}
We do not want to exclude  $k$-factor models with large $k$ a priori
if the data strongly indicates that such a complex model is appropriate.
Similarly, we do not want to exclude distributions with small $\lambda$.
We rather want to  avoid such models if the data does not give us strong enough
reason for  choosing them.
To solve this problem we use a generalized version of the structural
risk minimization principle (Vapnik 1998, Vapnik 1995).
In this way  the result of Theorem \ref{theo:Lern2} can be extended such that
$k$ and $\lambda$ need  not to be chosen in advanced. For doing so we have to 
define
a prior probability measure $\nu$  on the set
of possible pairs $(k,\lambda)$
in advance in order to get a bound on the true risk.

We obtain from Theorem \ref{theo:Lern2} with a standard union
bound argument, see e.g. Lemma 4.1 in  (Herbrich 2002),
following corollary.
\begin{corollary}
\label{theo:SRM}
  Let $(\lambda_n)_{n\in \N}$ be a strict monotone decreasing sequence
  converging to zero. Let $\nu$ be a strict positive probability measure on 
$\N^2$
  and $A_{kn} =\{1,\ldots k \} \times \{1,\ldots n\}$.
  Then for training data $\bx_1,\bx_2,\dots,\bx_l$ we have with probability at 
least $1-\eta$
\[
   R(P) \le R_{emp} (P) + \phi (k,\lambda_n, \eta \,\nu(A_{kn})) \,,
\]
for  $P \in {\cal H}_{k}^{\lambda_n}$ and $\phi$ defined above.
\end{corollary}
The probability measure $\nu$ in Corollary \ref{theo:SRM} has
to be chosen before seeing the
training data.  Usually one chooses $\nu$ to be inversely
proportional to $\phi$. However, it should be emphasized
that the statement in Corollary \ref{theo:SRM} does not assume
that the measure $\nu$ is chosen according to any prior
probabilities for the possible probability measures as used
in Bayesian learning.

Now we can formulate structural risk minimization as follows:
\[
\underset{(k,n)\in \N^2}{ \mbox{Minimize}}\quad  t  +  \phi
(k,\lambda, \eta \,\nu(A_{kn}))
\]
\[
\quad \mbox{subject to } \quad  t \;= \min_{\quad P \in {\cal 
H}_{k}^{\lambda_n}} R_{emp}(P)
\]

\subsection{Model Computation by Convex Optimization}
\label{sec:convex}
In this section we show that the minimization of the empirical risk  $R_{emp}$
is a convex optimization problem in the class of all $k$-factor models for a 
fixed $k$
and $\lambda$. We formulate the  optimization problem for $\ln P$ rather
than for $P$.
Let $\bx^1,\ldots,\bx^l$ be the training data.
For each set ${\bf j}=\{j_1,\dots,j_k\}\subset \{1,\dots,n\}$
let $V_{\bf j}$ be the vector space of real-valued functions
on the set
$
\Omega_{\bf j}:=\Omega_{j_1} \times \dots \times \Omega_{j_k}\,.
$
Define
\[
V^k:=\oplus_{\bf j} V_{{\bf j}}
\]
where ${\bf j}$ runs over all subsets of $\{1,\dots,n\}$ of size $k$.
Now fix $k$ and ${\bf j}$ of size $k$. The  set $\Omega_{\bf j}$ is discrete
and so the functions
\[
e_{w_1,w_2,\dots,w_k}(v_1,v_2,\dots,v_k) =
\begin{cases}
  1  & \text{if } \quad w_1 =v_1,\dots,w_k =v_k \\
  0  &   \text{else}
\end{cases}
\]
form a canonical basis of $V_{{\bf j}}$.
Further let ${\bf x}_{\bf j}$ be the $k$-tuple given by the restriction of ${\bf 
x} \in \Omega$
to the variables $x_{j_1},\dots,x_{j_k}$ and  ${\bf x}^i_{\bf j}$ the 
restriction of the training data
$\bx^i$ to the variables of ${\bf x}_{\bf j}$.

Now, let us formulate the optimization problem of the empirical risk.
Find a vector $f=\oplus_{\bf j} f_{\bf j} \in V^k$
with $f_{\bf j} \in V_{\bf j}$
such that
\begin{equation}
  \label{eq:GF}
   \tilde{R}_{emp}(f):= \frac{1}{l}\sum^l_{i=1} \sum_{\bf j} f_{\bf j} ({\bf 
x}^i_{\bf j})
\end{equation}
is minimal subject to the constraints
\[
Z(f):=\sum_{\bf x \in \Omega} \exp( \sum_{\bf j} f_{\bf j} ({\bf x}_{\bf j}))=1
\]
and
\[
\qquad \qquad G_{\bf x}(f):=\sum_{\bf j} f_l({\bf x}_{\bf j}) \geq \log \lambda
\quad \text{for all ${\bf x} \in \Omega$.}
\]
Then the probability measure $P \in \Hyp$ is given by $\ln f$.
The functions $\tilde{R}_{emp}$ and
$G_{\bf x}$ are obviously linear in $f$.
The function $Z$ is convex in $f$ since the exponent  is linear
in $f$ and the exponential function is convex. Clearly, the sum
of the convex functions
\[
f\mapsto \exp( \sum_{\bf j} f_{\bf j} ({\bf x}_{\bf j}))
\]
over all ${\bf x}$ is convex.
Therefore the problem can be solved with  standard convex optimization
methods (Pallaschke and Rolewicz 1997).
However, for a large number of variables the solution becomes computationally
expensive.
Notice that the number of constraints in (\ref{eq:GF}) increases
exponentially with $n$. Furthermore the calculation of
$\tilde{R}_{emp}(f)$ involves the summation over exponentially many
possible $n$-tuples ${\bf x}$.
This problem is similar to the calculation of the
partition function in statistical mechanics which requires a summation over
exponentially many states of a many-particle system.
Those computational problems are not  different
from those appearing in conventional  approaches
to log-linear model 
selection.%%%%%%%%%%%%%%%%%%%%%%%%%%%%%%%%%%%%%%%%%%%%%%%%%%%%%%%%%%%%%%%%%%%%%%
\section{CONCLUSION}

We used structural risk minimization of statistical learning theory to obtain
a new selection criterion for  log-linear models. In this we
way have not only a
criterion to choose the suitable model complexity, but
we found also a bound on the actual risk, i.e.
the key feature of our approach is that it
provides statements of the form ``if a log-linear model
of certain simplicity fits well to the observed data then
we are guaranteed with high probability that it will fit future observations
as well''.
This kind of nonasymptotic statistics  become more and more important
for in research areas like cognitive and social science, where the number of
possibly relevant features for the system under investigation is large.
It offers a new way of developing models, namely, infering statistical 
dependencies
amoung a large number of features from data.

We structured the class of log-linear models according to the degree of the 
interaction terms.
Using other measures of complexity than VC dimension like Gaussian
complexity (Bartlett and Mendelson 2001) we believe
that the obtained bounds can be improved. The structure of the model classes can 
be refined
taking for example prior knowledge into account.
In future this model selection criterion need to be tested on real-world data 
and the results
compared to existing ones.

%%%%%%%%%%%%%%%%%%%%%%%%%%%%%%%%%%%%%%%%%%%%%%%%%%%%%%%%%%%%%%%%%%%%%%%%%%%%%%%%
%%%%%%%%
%

\end{document}